\input amstex
\documentstyle{amsppt}
 \loadmsbm
  \loadbold

 \magnification=\magstep1
        \pagewidth{13cm}
        \pageheight{20cm}

\define\dQ{{\Bbb Q}}

\define\cF{{\Cal F}}
\define\cG{{\Cal G}}

\define\cO{{\Cal O}}

\define\cS{{\Cal S}}
\define\eq{{\frak q}}

\define\ep{{\frak p}}
\define\eP{{\frak P}}
\define\eQ{{\frak Q}}

\define\Gal{{\operatorname{Gal}}}

\define\res{{\operatorname{res}}}
\define\rank{{\operatorname{rank}}}
\define\Char{{\operatorname{char}}}

\define\tw#1{\,^{#1}\!}

\font\pfeile = cmsy10 scaled 1440
\newfam\pfeilfam
\textfont\pfeilfam=\pfeile
                   \scriptfont\pfeilfam=\pfeile
                                      \scriptscriptfont\pfeilfam=\pfeile

\mathchardef\swpfeil="2D2E

\scrollmode
\NoBlackBoxes

\define\dentzer{1}
\define\jz{2}
\define\kal{3}
\define\ks{4}
\define\neftin{5}
\define\nomura{6}
\define\plans{7}
\define\rw{8}
\define\schneps{9}
\define\serre{10}
\define\weirB{11}

\title
On the minimal ramification problem for $\ell$-groups
\endtitle
\author
Hershy Kisilevsky and Jack Sonn
\endauthor
\affil Concordia University, Montreal and Technion, Haifa
\endaffil
\address
Department of Mathematics, Concordia University, 1455 de Maisonneuve
Blvd. West, Montreal, Canada, and Department of Mathematics,
Technion, 32000 Haifa, Israel
\endaddress
\email kisilev\@mathstat.concordia.ca, sonn\@math.technion.ac.il
\endemail
\thanks
The research of the first author was supported by grants from
NSERC and FQRNT.
\endthanks
\thanks
The research of the second author was supported by the Fund for
the Promotion of Research at the Technion and the E. and M.
Mendelson Research Fund.
\endthanks


\abstract Let $\ell$ be a prime number.  It is not known if every
finite $\ell$-group of rank $n>1$ can be realized as a Galois
group over $\dQ$ with no more than $n$ ramified primes.  We prove
that this can be done for the (minimal) family of finite
$\ell$-groups which contains all the cyclic groups of $\ell$-power
order, and is closed under direct products, (regular) wreath
products, and rank-preserving homomorphic images. This family
contains the Sylow $\ell$-subgroups of the symmetric groups and of
the classical groups over finite fields of characteristic not
$\ell$.  On the other hand, it does not contain all finite
$\ell$-groups.
\endabstract

\endtopmatter
\document
\head  
Introduction \endhead 

Let $K$ be a global field, $L/K$ a finite Galois extension with Galois group $G=G(L/K)$.  Let $\ep$ be a finite
prime of $K$.  If $\ep$ ramifies in $L$ and if $\eP$ is a prime of $L$ dividing $\ep$, then the inertia group
$T(\eP/\ep)$ is a nontrivial subgroup of $G$.  If $T$ is the subgroup of $G$ generated by all $T(\eP/\ep)$, then
the fixed field of $T$ is an unramified extension of $K$.  If $K=\dQ$ then by Minkowski's theorem, there are no
nontrivial unramified algebraic extensions of $\dQ$, so $T=G$.  Suppose in addition that $L/\dQ$ is tamely
ramified, \it i.e. \rm for every prime $p$ ramified in $L/\dQ$,  all the $T(\eP/p)$ are cyclic of order prime to
$p$.  It follows in particular that if for each ramified $p$ we fix an inertia group $T(\eP/\ep)=\langle g_p
\rangle$, then the normal subgroup of $G$ generated by the $g_p$ is all of $G$.
\medskip We are interested in
the case where $G=G(L/\dQ)$ is an $\ell$-group, $\ell$ a prime.  Here $L/\dQ$ is
tamely ramified if and only if all the primes $p$ that ramify in $L$
are prime to $|G|$. Let $\bar G
=G/\Phi(G)$ be the quotient of $G$ by its Frattini subgroup
$\Phi(G)$.  Then the normal subgroup of $G$ generated by the $g_p$
is all of $G$ if and only if their images $\bar g_p$ in $\bar G$
generate $\bar G$ if and only if (by Burnside's basis theorem) the
$g_p$ generate $G$.  It follows that $\rank (G)$, the minimal number
of generators of $G$, is less than or equal to the number of primes
$p$ that ramify in $L$, or equivalently, the number of primes that
ramify in $L$ is at least $\rank(G)$.
\medskip It is an open
problem whether or not every finite $\ell$-group $G$ can be realized
as the Galois group of a tamely ramified extension of $\dQ$ with
exactly $\rank(G)$ ramified primes (see e.g. Plans \ \cite \plans).  We call this \it the minimal ramification problem.  \rm Using
Dirichlet's theorem on primes in arithmetic progressions, it is easy
to show that this problem has an affirmative answer for abelian
$\ell$-groups $G$.  It has been remarked in Serre \ \cite \serre \
that for odd $\ell$, the Scholz-Reichardt method for realizing
$\ell$-groups over $\dQ$ yields realizations of an $\ell$-group of
order $\ell^n$ with no more than $n$ ramified primes.  But
$n=\rank(G)$ only if $G$ is elementary abelian.  B. Plans \cite
\plans \  has improved this bound by showing that the
Scholz-Reichardt method yields the bound equal to the sum of the
ranks of the factors of the lower central series of $G$ (without
the bottom factor). Thus the minimal ramification problem has an affirmative solution for
odd order $\ell$-groups $G$ of nilpotency class $2$.  A. Nomura
\cite \nomura \  has refined Plans' result and proved that the minimal ramification
problem has an affirmative solution for $3$-groups of order $\leq
3^5$.

\medskip In this paper we produce (for every $\ell$,
including $\ell=2$) a new family of $\ell$-groups for which the minimal ramification problem has an affirmative
solution. To be precise, given the prime $\ell$, let $\cG(\ell)$ be the minimal family of $\ell$-groups containing
the cyclic $\ell$-groups and which is closed under direct products, (regular) wreath products, and rank-preserving
homomorphic images. Then every group $G$ in $\cG(\ell)$ is tamely realizable over $\dQ$ with exactly $\rank(G)$
ramified (finite) primes. The family $\cG(\ell)$ contains all direct products of iterated wreath products of
cyclic groups of $\ell$-power order, and in particular, all Sylow $\ell$-subgroups of the symmetric groups (Kaloujnine \ \cite
\kal) \ and of the classical groups over finite fields of characteristic prime to $\ell$ (Weir \ \cite \weirB).  On the
other hand, it does not contain all finite $\ell$-groups, as we will see.
\head  
$\ell$-groups as Galois groups with minimal ramification \endhead
\medskip Let $G,H$ be finite (abstract) groups.  We will define the (regular) \it wreath product \rm
$H\wr G$ of $H$ with $G$ to be the semidirect product $H^{|G|}
\rtimes G$, where $H^{|G|}$ is the direct product of $|G|$ copies
of $H$, with $G$ acting on $H^{|G|}$ by permuting the copies of
$H$ like the regular (Cayley) representation of $G$. Define the
$n^{\text {th}}$ iterated wreath product $G^{\wr n}$ of $G$ by
$G^{\wr 1}:=G$ and $G^{\wr n}:=  G^{\wr (n-1)}\wr G$ for $n>1$.

\bigskip
\bf Proposition 1. \rm (Ribes-Wong \ \cite \rw ) \ \ \it Let $G,H$ be finite $\ell$-groups of ranks
$m,n$ respectively.  Then  $\rank(H \wr G)=m+n$.
\smallskip Proof. \rm Let $G$ have minimal generating set  $\{
g_1,...,g_m \}$ and let $H$ have minimal generating set
$\{h_1,...,h_n \}$.  Then it is clear that $H \wr G$ is generated
by  $\{ g_1,...,g_m,h_1,...,h_n \}$, so $\rank(H \wr G) \leq m+n$.
Now if $\rank(H \wr G)<m+n$, then by Burnside's basis theorem, a
proper subset of $\{ g_1,...,g_m,h_1,...,h_n \}$ would generate $H
\wr G$.  But if a $g_i$ is dropped from this generating set, the
resulting subgroup is of the form $H \wr G_1$ with $G_1$ a proper
subgroup of $G$, so $H \wr G_1$ is a proper subgroup of $H \wr G$.
Similarly, if an $h_i$ is dropped from this generating set, the
resulting subgroup is of the form $H_1 \wr G$ with $H_1$ a proper
subgroup of $H$, so $H_1 \wr G$ is a proper subgroup of $H \wr G$.
\qed
\bigskip

We will say that an extension of global fields $L/K$ contains no
non-trivial unramified subextension, or that $L$ contains no
non-trivial unramified subextension of $K,$  if whenever $K \subseteq E \subseteq L$
are field extensions with $E/K$ unramified then $E=K.$
\bigskip

Fix an arbitrary global field $k$ and a prime $\ell\neq\text
{char}(k)$. Define a family $\cF^{\min}:=\cF^{\min}_{k,\ell}$ of
(isomorphism classes of) finite $\ell$-groups as follows:
$G\in\cF^{\min}$ iff given any finite set $S$ of primes of $k$ and
any finite separable extension $K/k$, there exists a finite Galois
extension $L/K$ with $G(L/K)\cong G$ such that the set of primes
$\{\ep_1,...,\ep_n\}$ of $K$ that ramify in $L$ satisfy the
following $5$ conditions:
 \smallskip

1.\ \ $n=\rank(G)$ (the minimal number of generators of $G$)
\smallskip

2.\ \ the primes $p_1,...,p_n$ of $k$ below $\{\ep_1,...,\ep_n\}$
are distinct \smallskip

3.\ \ $\{p_1,...,p_n\}\cap S =\emptyset $
\smallskip

4.\ \ $p_1,...,p_n$ split completely in $K$ \smallskip

5.\ \ $L$ contains no non-trivial unramified subextension of $K.$
\smallskip

The main result of this paper is the following:

\bigskip

{\bf Theorem 1.}  {\it The family $\cF^{\min}$ has the following
properties:
\smallskip

a) \ \ $\cF^{\min}$ contains all cyclic groups of $\ell$-power order.
\smallskip

b) \ \ If $G,H \in \cF^{\min}$, then $G\times H \in \cF^{\min}$.
\smallskip
c) \ \ If $G \in \cF^{\min}$ and $N$ is a normal subgroup of $G$ contained
in the Frattini subgroup $\Phi(G)$ of $G$, then $G/N \in \cF^{\min}$.
\smallskip
d) \ \ If $G,H \in \cF^{\min}$, then $H \wr G \in \cF^{\min}$.}

\bigskip

Before proving the theorem we note the following immediate
consequence when $k=K=\dQ$.

\bigskip
{\bf Corollary 1.} {\it Let $\cG(\ell)$ be the minimal family of
$\ell$-groups satisfying conditions a)-d) of Theorem 1, i.e.
$\cG(\ell)$ contains all cyclic groups of $\ell$-power order and is
closed under direct products, (regular) wreath products, and
rank-preserving homomorphic images. Then all $G\in\cG(\ell)$ of rank
$n$ are tamely realizable over $\dQ$ with exactly $n$ ramified
primes.}

\bigskip

We will use the following lemma in the proof of Theorem 1.

\bigskip

\proclaim{Lemma 1} Suppose that $K_1$ and $K_2$ are Galois
extensions of $K$ with $\Gal(K_i/K)=G_i, i=1,2,$ and such that
$K_2/K$ contains no non-trivial unramified subextensions. Suppose
also that the extensions $K_1/K$ and $K_2/K$ are ramified at
disjoint sets of primes of $K.$ Then $K_1\cap K_2=K,$ \rm(hence
$G=\Gal(K_1\cdot K_2/K)\cong G_1\times G_2$), \it and for any
unramified subextension  $K \subseteq E \subseteq K_1\cdot K_2,$ we
have $K \subseteq E \subseteq K_1.$ In particular if $K_1/K$ also
contains no non-trivial unramified subextensions then $K_1\cdot
K_2/K$ contains no non-trivial unramified subextensions.

\endproclaim

\bigskip
{\it Proof.} Let $\{\ep_1,\ldots,\ep_s\}$ be the primes of $K$
ramified in $K_1$ and  $\{\eq_1,\ldots,\eq_t\}$ be the primes of
$K$ ramified in $K_2.$ Then by assumption
$\{\ep_1,\ldots,\ep_s\}\cap\{\eq_1,\ldots,\eq_t\}=\emptyset .$
Since $K_1\cap K_2 \subseteq K_1,$  we see that $K_1 \cap K_2/K$
is ramified only at primes in $\{\ep_1,\ldots,\ep_s\},$ and
similarly $K_1\cap K_2 \subseteq K_2,$ implies that $K_1 \cap
K_2/K$ is ramified only at primes in $\{\eq_1,\ldots,\eq_t\}.$
Therefore $K_1 \cap K_2/K$ is unramified, and since $K_2/K$
contains no non-trivial unramified subextension we see that $K_1
\cap K_2 =K,$ and so $\Gal(K_1\cdot K_2/K)\cong G_1\times G_2.$
Let $T_{\eQ}\subseteq G=\Gal(K_1\cdot K_2/K)$ be the subgroup
generated by all the inertia groups $T(\eQ_i/\eq_i)$ where $\eQ_i$
runs over all primes of $K_1\cdot K_2$ dividing some prime $\eq_i
\in \{\eq_1,\ldots,\eq_t\}.$ Since $K_1/K$ is unramified at the
primes $\{\eq_1,\ldots ,\eq_t\},$ we see that $K\subseteq
K_1\subseteq (K_1\cdot K_2)^{T_{\eQ}}.$ But since $G\cong
G_1\times G_2,$ we have that the restriction map $\res:\Gal
(K_1\cdot K_2/K_1)\longrightarrow G_2$ is an isomorphism. Also
since $K_2/K$ contains no non-trivial unramified subextension, it
follows that $\res(T_{\eQ})=G_2,$ and therefore
$T_{\eQ}=\Gal(K_1\cdot K_2/K_1)$ and $ K_1= (K_1\cdot
K_2)^{T_{\eQ}}.$ Suppose that $K \subseteq E \subseteq K_1\cdot
K_2$ with $E/K$ unramified. Then $E$ is contained in the subfield
of $K_1\cdot K_2$ fixed by $T_{\eQ}$. But then $E$ is fixed by
$T_{\eQ}$ and therefore $E\subseteq K_1.$ If $K_1/K$ contains no
non-trivial unramified subextension, we must have $E=K.$ \qed


\bigskip
\par\noindent
We will also need a lemma from Kisilevsky-Sonn \cite \ks:
\bigskip
\par\noindent
Let $K$ be a global field, $\frak p$ a finite prime of $K$,
$I_{\frak p}$ the group of fractional ideals prime to $\frak p$,
$P_{\frak p}$ the group of principal fractional ideals in $I_{\frak
p}$, $P_{\frak p,1}$ the group of principal fractional ideals in
$P_{\frak p}$ generated by elements congruent to $1$ mod $\frak p$.
Then $Cl_K = I_{\frak p}/P_{\frak p}$ is the class group of $K$,
$Cl_{K,\frak p}=I_{\frak p}/P_{\frak p,1}$ is the ray class group
with conductor $\frak p$, and $\overline P_{\frak p}=P_{\frak
p}/P_{\frak p,1}$ is the principal ray with conductor $\frak p$. We
have a short
 exact sequence
$$ 1 \longrightarrow \overline P_{\frak p} \longrightarrow Cl_{K,\frak p}
\longrightarrow Cl_K \longrightarrow 1. \tag*$$
\bigskip
\par\noindent
For prime $\ell\neq \Char(K)$ we consider the exact sequence of
$\ell$-primary components
$$ 1 \longrightarrow \overline P_{\frak p}^{(\ell)} \longrightarrow
Cl_{K,\frak p}^{(\ell)} \longrightarrow Cl_K^{(\ell)} \longrightarrow 1.
\tag*{${}_{\ell}$}$$
 \bigskip
\par\noindent
We are interested in primes $\frak p$ for which the sequence (*${}_{\ell}$)
splits. Let $\frak a_1,...,\frak a_s \in I_K$ be such that their images
$\overline{\frak a}_i$ in $Cl_K^{(\ell)}$ form a basis of the finite abelian
$\ell$-group $Cl_K^{(\ell)}$.  Let $\ell^{m_i}$ be the order of
$\overline{\frak a}_i$, $i=1,...,s$.  Then
$\frak a_i^{\ell^{m_i}}=(a_i)\in P_K$,
$i=1,...,s$.  Let
$K'=K(\zeta_{\ell^m},\root{\ell^{m}}\of{\epsilon},\root{\ell^{m_i}}\of{a_i},
1\leq i \leq s),$ the field extension obtained by adjoining
$\zeta_{\ell^m}$ a primitive $\ell^m$th root of unity,
the $\ell^m$th roots of all units $\epsilon$ of $K$, and the
$\ell^{m_i}$th roots of the elements $a_i\in K,$
where $m\geq \max\{1,m_1,...,m_s\}$.
\bigskip
\par\noindent
\proclaim {Lemma 2} {\rm (``Splitting Lemma") (Kisilevsky-Sonn \cite {\ks, Lemma
2.1})} In order that the sequence (*${}_{\ell}$) split, it is
sufficient that $\frak p$ split completely in $K'$.
\endproclaim

(For the proof see \cite \ks.)

\bigskip

{\bf Proof of Theorem 1.}

\bigskip

Let $K,S$ be given.

\smallskip
a) \ \ Let $p\notin S$ be a prime of $k$ which splits completely in $K'$, where
$K'$ is the field defined in the Splitting Lemma for $K.$
Let $\ep$ be a prime of $K$ dividing $p$.
Then by the Splitting Lemma the $\ell$-ray class field $R_{\frak p}$
of $K$ belonging
to the ray class group $Cl_{K,\frak p}^{(\ell)}$ has Galois group isomorphic to
$ Cl_K^{(\ell)}\times  \overline P_{\frak p}^{(\ell)}.$
Since the ${\ell}$-Hilbert Class field $H_K^{(\ell)}$ belongs
to $ Cl_K^{(\ell)}$
we see that
$R_{\frak p}=H_K^{(\ell)}\cdot L'$ with $H_K^{(\ell)}\cap L'=K,$ and
$\Gal(L'/K)\cong  \overline P_{\frak p}^{(\ell)}.$
Under our assumption that all units are $\ell^m$th powers
modulo $\ep,$ it follows that
$$
 \overline P_{\frak p}^{(\ell)}/ (\overline P_{\frak p}^{(\ell)})^{\ell^m}
\cong (\cO_K/\ep)^{\ast}/ ((\cO_K/\ep)^{\ast})^{\ell^m}
$$
is cyclic and has order divisible by $\ell^m.$ Taking $m \geq r,$
we see that there exists a cyclic extension $L/K$ of degree
$\ell^r$ ramified only at $\ep$, and in which $\ep$ is totally
ramified. $L/K$ satisfies conditions 1.-5. (with $n=1$).
 \smallskip

b)\ \  Since $G \in \cF^{\min},$ there is an extension $K_1/K$
with $\Gal(K_1/K) \cong G$ and which satisfies properties 1.-5.
with primes $\{p_1,...,p_n\},$ and  $\{\ep_1,...,\ep_n\}$
respectively. Set  $S'=S\cup\{ p_1,...,p_n\}.$ Since $H \in
\cF^{\min},$ let $K_2/K$ be an extension with $\Gal(K_2/K) \cong
H$ and which satisfies properties 1.-5. for $K,S',$ with primes
$\{q_1,...,q_m\},$ and  $\{\eq_1,...,\eq_m\}$ respectively. Then
by Lemma 1, $L=K_1K_2$ puts $G\times H$ in $\cF^{\min}$ with
primes $\ep_1,...,\ep_n, \eq_1,...,\eq_m,$ and where
$n+m=\rank(G\times H)$. This proves b).
\smallskip

c)\ \ Let $L/K$ be a Galois extension with group $G$ which puts
$G$ in $\cF^{\min}$.  Let $N$ be a normal subgroup of $G$
contained in $\Phi(G)$.  Let $L'$ be the fixed field of $N$. Then
$\rank(G/N)=\rank(G)$.  The other conditions are immediate.
\smallskip

d)\ \  Let $K_1/K$ be a Galois extension with group $G$ which puts
$G$ in $\cF^{\min}$, with ramified primes $\ep_1,...,\ep_n$ over
$p_1,...,p_n \notin S$.  Let $m=\rank(H)$, $S_1=S \cup \{
p_1,...,p_n\}$.  Apply the hypothesis $H\in \cF^{\min}$ to the
pair $K_1,S_1$.  Then there exists a Galois extension $L_1/K_1$
with group $H$, with $m$  primes $\eQ_{1},...,\eQ_{m}$ of $K_1$
ramified in $L_1$ such that the primes $q_{1},...,q_{m}$ of $k$
below $\eQ_{1},...,\eQ_{m}$ are distinct, $q_{1},...,q_{m}$ split
completely in $K_1$, $q_{1},...,q_{m} \notin S_1$, and $L_1$
contains no non-trivial unramified extension of $K_1$.  Let
$\eq_{1},...,\eq_{m}$ be the primes of $K$ below
$\eQ_{1},...,\eQ_{m}$. Then $\eq_{1},...,\eq_{m}$ split completely
in $K_1$.  So each $\eQ_i$ has $|G|$ distinct conjugates
$\{\sigma(\eQ_i)|\sigma \in G\}$ over $K$, $i=1,...,m$.  For each
$\sigma \in G$ the conjugate extension $\sigma(L_1)/K_1$ is
well-defined, since $L_1/K_1$ is Galois.  Let $L$ be the composite
of the $\sigma(L_1), \sigma \in G$.  For each $\sigma \in G$,
$\sigma(L_1)/K_1$ is Galois with group $H$, with exactly $m$
ramified primes $\sigma(\eQ_{1}),...,\sigma(\eQ_{m})$ lying above
$q_{1},...,q_{m}$, and $\sigma(L_1)$ contains no unramified
extension of $K_1$. Furthermore, the set of primes
 $\sigma(\eQ_{1}),...,\sigma(\eQ_{m})$ ramified in  $\sigma(L_1)/K_1$
is disjoint from the set of primes  $\tau(\eQ_{1}),...,\tau(\eQ_{m})$
ramified in $\tau(L_1)/K_1$ if $\sigma \neq \tau.$ This follows since
if $\sigma(\eQ_{i})=\tau(\eQ_{j})$ we would have $q_i=q_j.$ But then
$i=j$ by property 3 in the definition of $\cF^{\min}$ and so would have $\sigma=\tau.$
  \smallskip Applying Lemma 1
repeatedly we see that the fields $\{\sigma(L_1):\sigma \in G \}$
are linearly disjoint over $K_1$. It follows that we have an exact
sequence of groups $$ {\text(*)} \ \ \ 1 \rightarrow H^{|G|}
\rightarrow G(L/K) \rightarrow G \rightarrow 1 $$  where $G$ is
identified with $G(K_1/K)$ and $H^{|G|}$ is the direct product of
$|G|$ copies of $H$. Furthermore, this exact sequence defines a
unique homomorphism (injective in this case) $\phi:G\rightarrow
Out(H^{|G|})$, which is equivalent, as a permutation
representation on the $|G|$ copies of $H$, to the regular
representation of $G$. The set of all group extensions of $G$ by
$H^{|G|}$ corresponding to a given $\phi$, if nonempty, is in
one-one correspondence with $H^2(G,Z(H^{|G|}))$, (see Johnson-Zassenhaus \ \cite
\jz), where $Z(H^{|G|})$ denotes the center of $H^{|G|}$.  Since
$Z(H^{|G|})=Z(H)^{|G|}$ is an induced $G$-module,
$H^2(G,Z(H^{|G|}))=0$.  It follows that the group extension (*)
splits, and $G(L/K)\cong H \wr G$.

\smallskip The primes of $K$ that
ramify in $L$ are exactly
$\{\ep_1,\ldots,\ep_n,\eq_1,\ldots,\eq_m\}$, $n+m=\rank(H\wr G)$,
the primes $p_1,...,p_n,q_1,\ldots,q_m$ below
$\ep_1,\ldots,\ep_n,\eq_1,\ldots,\eq_m,$ are distinct, split
completely in $K$, and lie outside $S$.  Finally, $L/K$ does not
contain a non-trivial unramified subextension  $M/K$, since if it
did, then $M$ would be contained in $K_1$ and $K_1/K$ contains no
non-trivial unramified subextension of $K$. \qed
\bigskip
How large is the family $\cG(\ell)$?  It is smaller than the
family of all $\ell$-groups, as we will now show.

\bigskip \bf
Lemma 3. \it Let $G$ be a nontrivial group in $\cG(\ell)$, and let $dl(G)$  be the derived length
\rm(length of the derived series) \it of $G$. Then $dl(G)\leq \rank(G)$.
\medskip

Proof.  \rm By induction on the minimal number $t$ of applications of the three types of operations
(direct product, wreath product, rank-preserving homomorphic image) defining $\cG(\ell)$, which are
needed to produce $G$ starting from cyclic $\ell$-groups.  If $t=0$ ($G$ cyclic), we have $
dl(G)=\rank(G)$.   We examine the behavior of the rank and of the derived length under each of the
three operations.
\smallskip (i) If $G,H\in \cG(\ell)$, then $\rank(G\times
H)=\rank(G)+\rank(H)$ while $dl(G\times H)=\max(dl(G),dl(H))$
\smallskip (ii) $G,H\in \cG(\ell)$, then $\rank(H\wr
G)=\rank(G)+\rank(H)$ (Proposition 1) while $dl(H\wr G)\leq
dl(G)+dl(H))$ (easy)
\smallskip (iii) If $G \in \cG(\ell)$ and $\overline G$ is a homomorphic
image of $G$ (with $\rank(\overline G)=\rank(G)$), then $dl(\overline G)\leq dl(G)$.
The result follows.  \qed
\bigskip
\bf Proposition 2. \it For every $\ell$ and $n>1$, there exist
$\ell$-groups of rank $n$ not in $\cG(\ell)$.
\medskip
Proof.  \rm It suffices to see that for every $n>1$, there exist
$\ell$-groups of rank $n$ and derived length larger than $n$.  Let
$F$ be the free group of rank $n$, and let $F_t$ be the $t$th term
of the descending $\ell$-central series of $F$ ($F_1=F$ and for
$t>1$, $F_t=F_{t-1}^{\ell}[F,F_{t-1}]$).  It suffices to show that
the derived length of $F/F_t$ is larger than $n$ for sufficiently
large $t$.  But this is true since the derived length of $F$ is
infinite and the descending $\ell$-central series of $F$ has trivial
intersection.  (For sufficiently large $t$, $F_t$ does not contain
the (nontrivial) $n$th term of the derived series of $F$.) \qed

\bigskip \it Example 1.  \rm Here is an example of an $\ell$-group not in the family
$\cG(\ell)$. (We thank John Labute for help with this example.)
\smallskip Let $F$ be a free group on 2 generators $x,y$ and let $G$
be the quotient of $F$ by the sixth term $F_6$ of the descending
$\ell$-central series of $F$.  We claim $G \notin \cG(\ell)$. By
Lemma 3, it suffices to show that $dl(G)=3$. Indeed,
$[[x,y],[x,[x,y]]]$ lies in $F_5$ but not in $F_6$, so there are
two elements of the commutator subgroup $G'$ of $G$ whose
commutator is nontrivial.  (For another example see Remark 2
below.)
\bigskip
\it Remark 1. \rm If we drop condition 1. from the definition of
$\cF^{\min}$ to obtain the (larger) family $\cF$, then we get the
following variant of Theorem 1:

\medskip
{\bf Theorem 2.}  {\it The family $\cF$ has the following
properties:
\smallskip

a) \ \ $\cF$ contains all cyclic groups of $\ell$-power order.
\smallskip

b) \ \ If $G,H \in \cF$, then $G\times H \in \cF$.
\smallskip
c) \ \ If $G \in \cF$, then every homomorphic image of $G$ is in $
\cF$.
\smallskip
d) \ \ If $G,H \in \cF$, then $H\wr G \in \cF$.}
\medskip
The proof is the same as that of Theorem 1, \it mutatis mutandis.
\rm As with Theorem 1, we obtain

\medskip
\bf Corollary 2. \it Let $\hat{\cG}(\ell)$ be the minimal family
of $\ell$-groups satisfying conditions a)-d) of Theorem 2.   Then
all $G\in\hat{\cG}(\ell)$ are tamely realizable over $\dQ$.
\bigskip \rm Theorem 2 in fact gives tame realizations of the
groups in $\hat{\cG}(\ell)$ over every global field, which of
course follows from the Scholz-Reichardt theorem for $\ell$ odd,
and from Shafarevich's theorem for $\ell=2$. However for these
groups we obtain a different, perhaps simpler, proof, especially
for $\ell=2$.
\medskip
\it Remark 2.  \rm A finite group $G$ is called \it semiabelian \rm
 if and only if there exists a sequence
$$G_0=\{1\},G_1,...,G_n=G$$ such that $G_i$ is a homomorphic image
of a semidirect product $A_i\rtimes G_{i-1}$ with $A_i$ abelian, $i=1,...,n$. \medskip\rm  It turns
out that $\hat{\cG}(\ell)$ is the family of all  semiabelian \rm $\ell$-groups, as
we will show.   Dentzer \ \cite \dentzer \ gives geometric realizations of the semiabelian groups
over $k(t)$ for any field $k$, in particular for $k$ a global field, so by Hilbert's irreducibility
theorem, realizations over global fields $k$. However it does not seem to be known how to produce
tame realizations via Hilbert's irreducibility theorem. Dentzer \ \cite \dentzer \ also gives an
example of a three generator $\ell$-group of order $\ell^5$ (for any odd $\ell$) which is not
semiabelian.
\smallskip  \bf Proposition 3. \it For any prime $\ell$,
$\hat{\cG}(\ell)$ is the family of all semiabelian $\ell$-groups.
\medskip Proof.  \rm Let $\cS(\ell)$ denote the family of all
semiabelian $\ell$-groups.  It is clear from the definition that $\cS(\ell)$ contains all cyclic
$\ell$-groups and is closed under homomorphic images.  Furthermore by Dentzer \ \cite {\dentzer, Theorem
2.8}, $\cS(\ell)$ is closed under direct products and (regular) wreath products.  Hence $\cS(\ell)$
contains $\hat{\cG}(\ell)$.  For the reverse inclusion, suppose contrarily that $G$ is a group of
minimal order in $\cS(\ell)\setminus \hat{\cG}(\ell)$.  Then $G$ is nonabelian, hence nontrivial.
By Dentzer \ \cite {\dentzer, Theorem 2.3}, $G$ is a composite $AH$ with $H$ a proper semiabelian subgroup of
$G$ and $A$ an abelian normal subgroup of $G$.  Then $G$ is a homomorphic image of a semidirect
product $A\rtimes H$, and by the induction hypothesis, $H \in \hat{\cG}(\ell)$.  Now $A\rtimes H$
is a homomorphic image of the (regular) wreath product $A\wr H$, which lies in $\hat{\cG}(\ell)$,
and hence so does its homomorphic image $AH=G$, a contradiction.  \qed \medskip
\bigskip
\it Remark 3.  \rm Given a finite $\ell$-group $G$, let $ram^t(G)$ denote the minimal $n$ such that
$G$ can be realized as a Galois group of a tamely ramified extension $L/\dQ$ with exactly $n$
ramified primes. As mentioned in the introduction, Plans \ \cite \plans \ has shown that the
Scholz-Reichardt method for realizing odd order $\ell$-groups over $\dQ$ can be made to yield an
upper bound for $ram^t(G)$ equal to the sum of the ranks of the factors in the lower central series
of $G$, where the bottom factor can be left out of this sum.  For most of the groups in the family
$\cG(\ell)$, this bound is larger than the rank of the group, e.g. for $C_{\ell}\wr C_{\ell}$,
$\ell>3$. \bigskip
 \bf Note.  \rm Since the submission of this paper, it has been proved by D.
Neftin \cite \neftin \ that the family ${\cG}(\ell)$ is \it equal \rm to the family
$\hat{\cG}(\ell)$ of semiabelian $\ell$-groups. To give some indication of the size of
${\cG}(\ell)$, the following is known  about ``small" $\ell$-groups (Dentzer \cite \dentzer; see also Schneps \ \cite
\schneps):
\smallskip
1. \ \ For any $\ell$, all $\ell$-groups of order $\leq \ell^4$ are semiabelian.
\smallskip
2. \ \ All $2$-groups of order $\leq 32$ are semiabelian. \smallskip 3. \ \ Among the $267$ groups
of order $64$, only $10$ are not semiabelian.  Similarly,  among the $2328$ groups of order $2^7$,
 $82$ are not semiabelian, and among the $56092$ groups of order $2^8$,  $993$ are not
semiabelian. Among the $67$ groups of order $3^5$, $10$ are not semiabelian, and among the $504$
groups of order $3^6$, $54$ are not semiabelian.



\bigskip \bigskip

 \Refs

 \ref \key \dentzer \by R. Dentzer \paper On geometric embedding
 problems and semiabelian groups \jour Manuscripta Math. \vol 86
 \yr 1995 \pages 199-216 \endref


\ref \key \jz \by C. E. Johnson, and H.  Zassenhaus
\paper On equivalence of finite group extensions \jour Math. Z.
\vol 123 \yr 1971 \pages 191--200
\endref

\ref \key \kal \by L. Kaloujnine \paper La structure des
 $p$-groupes de Sylow des groupes symetriques finis \jour Ann. ecole
 Norm. \vol 65 \yr 1948 \pages 239-276 \endref



\ref \key \ks
\by H. Kisilevsky and J. Sonn
\paper Abelian extensions of global fields with constant local degrees
\jour Math. Research Letters
\vol 13 {\rm{no. 4}}
\yr 2006
\pages 599--605
\endref

\ref \key \neftin \by D. Neftin \paper On semiabelian $p$-groups \rm (preprint) \endref

\ref \key \nomura \by A. Nomura \paper Notes on the minimal number
of ramified primes in some $l$-extensions of $\dQ$ \jour Archiv.
Math. \vol 90 \yr 2008 \pages 501-510 \endref

 \ref \key \plans \by B. Plans \paper On the minimal number of
 ramified primes in some solvable extensions of $\dQ$  \jour Pacific
 J. Math. \vol 215  \yr 2004 \pages 381-391 \endref

\ref \key \rw \by L. Ribes and K. Wong \paper On the minimal number of generators of
certain groups, \rm in  Groups St. Andrews 1989 \jour LMS Lecture Notes \endref


\ref \key \schneps \by L. Schneps \paper Reduction of $p$-groups \jour Comm. Alg. \vol 21  \yr
1993 \pages 1603-1609 \endref

\ref \key \serre \by J.-P. Serre \book Topics in Galois Theory
\publ Jones and Bartlett \publaddr Boston \yr 1992 \endref

 

\ref \key \weirB \by A. J. Weir \paper Sylow $p$-subgroups of the
classical groups over finite fields with characteristic prime to $p$
\jour Proc. Amer. Math. Soc. \vol 6 \yr 1955 \pages 529-533 \endref

\endRefs
\enddocument